\documentclass[runningheads]{llncs}

\usepackage{amsmath,amsfonts,amssymb,amsopn}

\usepackage[T1]{fontenc}
\usepackage{graphicx}
\usepackage{color}
\usepackage{stmaryrd}

\newcommand{\mybrL}{\Big\llbracket}
\newcommand{\mybrR}{\Big\rrbracket}

\renewcommand*{\Omega}{\varOmega} 
\newcommand{\sub}{\scriptscriptstyle} 
\renewcommand*{\vec}{\boldsymbol} 
\newcommand{\df}{\mathrm{d}} 
\newcommand{\dt}{\df t} 
\newcommand{\sdbt}{\vec{\sigma} \df \vec{B}_t} 

\newcommand{\grad}{\boldsymbol{\nabla}} 

\newcommand{\tp}{^{\scriptscriptstyle T}} 
\renewcommand*{\div}{\boldsymbol{\nabla\cdot}} 

\newcommand{\adv}{\boldsymbol{\cdot\nabla}} 

\newcommand{\alf}{\frac{1}{2}} 

\newcommand{\Den}{\operatorname{Den}}
\newcommand{\M}{\mathbb{M}}
\newcommand{\R}{\mathbb{R}}

\newcommand{\Div}{{\rm Div}}

\newcommand{\C}{\mathcal{C}}

\begin{document}
\title{Casimir-dissipation stabilized stochastic rotating shallow water equations on the sphere}
\titlerunning{Casimir-dissipation stabilized stochastic RSW equations}
%
\author{Werner Bauer\inst{1}\orcidID{0000-0002-5040-4287} \and
R\"udiger Brecht\inst{2}\orcidID{0000-0003-4040-3122} 
}
\authorrunning{W. Bauer and R. Brecht}
%
\institute{University of Surrey, Guildford, GU2 7XH, UK \\
\email{w.bauer@surrey.ac.uk}\\
\and
Universit\"at Hamburg, Hamburg, Germany
}
\maketitle              
\begin{abstract}
We introduce a structure preserving discretization of stochastic 
rotating shallow water equations, stabilized with an energy conserving
Casimir (i.e. potential enstrophy) dissipation. A stabilization of a stochastic scheme is usually required as, by modeling subgrid effects
via stochastic processes, small scale features are injected which often lead to noise on the grid scale and numerical instability. Such noise is usually 
dissipated with a standard diffusion via a Laplacian which necessarily also dissipates energy.
In this contribution we study the effects of using an energy preserving selective Casimir dissipation method compared to diffusion via a Laplacian. For both, we analyze stability and accuracy of the stochastic scheme. The results for a test case of a barotropically unstable jet show that Casimir dissipation allows for stable simulations that preserve energy and exhibit more dynamics than comparable runs that use a Laplacian. 

\end{abstract}

\section{Introduction}

Stochastic modeling allows one to not only simulate the time evolution of a dynamical system but also to estimate the reliability of the prediction through ensemble runs. When estimating the related uncertainty, it takes into account errors in measurements, in the mathematical model, and via numerical approximations. This is an important advantage over traditional approaches of deterministic models where such uncertainty estimates are not available. 

Apart from uncertainty estimates, it is important to preserve quantities such as mass, energy, or circulation to guarantee that the statistics of the solutions is correctly represented and that the numerical schemes remain stable and accurate. Deriving consistent stochastic models that preserve such invariants is challenging, but can be achieved using e.g. the Location Uncertainty (LU) framework \cite{Memin2014}. For example, in \cite{Brecht2021b} we derived energy preserving stochastic rotating shallow water (RSW) equations. However, there are no satisfying structure preserving discretization methods available for such stochastic flow models.

As a first step towards a consistent, energy preserving discretization of these stochastic RSW equations, we applied in \cite{Brecht2021b} a variational discretization for the deterministic part (introduced in \cite{BauerRSW2019}) while treating the stochastic terms using standard finite difference (FD) operators. In the following, we refer to this discrete system as stochastic RSW-LU scheme. As shown in \cite{Brecht2021b}, for homogeneous noise where no additional diffusion (in form of a Laplace operator) was required, this system preserved spatially the energy (currently there is no suitable energy preserving stochastic time integrator available). However, when modeling subgrid effects via stochastic processes, small scale features are injected on the grid level leading to numerical instabilities, especially for inhomogeneous noise. Thus, we applied a biharmonic Laplacian for stabilization, which dissipates energy. 

As a remedy to this latter point, we suggest here to use Casimir dissipation for stabilization to avoid the dissipative effect of the Laplacian on the total energy while keeping the stochastic scheme stable also in case of inhomogenous noise. The suggested combination is possible because the derivation of the Casimir dissipation for the RSW equations in \cite{Brecht2021a} follows the same variational discretization principle of \cite{BauerRSW2019} applied to the RSW equations that form the deterministic core of the stochastic RSW-LU model of \cite{Brecht2021b}. 

In Section~\ref{sec_lu}, we introduce briefly the main idea and some steps of the derivation of the stochastic flow model. In Section~\ref{sec_discretization}, we explain shortly how the variational principle can be used to derive a scale selective Casimir dissipation, and we combine it with the stochastic model. In Section~\ref{sec_numerics}, we test the novel scheme on a standard test cases focusing on the energy preserving properties. In Section~\ref{sec_conclusion}, we draw some conclusions.

\section{Stochastic RSW equations under Location Uncertainty}
\label{sec_lu}

We briefly introduce the stochastic RSW equations under Location Uncertainty \cite{Memin2014}
(RSW-LU) used for this study. The full derivation can be found in \cite{Brecht2021b}.

Let $\vec{u}=(u,v)$ be the 2D velocity and $h$ the water depth of the RSW flow. The unresolved random flow is modeled by the random flow term $\sdbt$ (of null ensemble mean) that can be computed with the spectral decomposition: 
$	\vec{\sigma} (\vec{x}, t)\, \df \vec{B}_t = \sum_{n \in \mathbb{N}} \vec{\Phi}_n (\vec{x}, t)\, {\df \beta_t^n},\   \vec{a} (\vec{x}, t) = \sum_{n \in \mathbb{N}} \vec{\Phi}_n (\vec{x}, t) \vec{\Phi}_n\tp (\vec{x}, t),$
where $\beta^n$ denotes $n$ independent and identically distributed (i.i.d.) one-dimensional standard Brownian motions. 
Here, we used the fact that the noise $\sdbt$ and its variance $\vec{a}$ (measuring the strength of the noise) can be represented by an orthogonal eigenfunction basis $\{  \vec{\Phi}_n (\bullet, t) \}_{\sub n \in \mathbb{N}}$ weighted by the eigenvalues $\Lambda_n \geq 0$ such that $\sum_{\sub n \in \mathbb{N}} \Lambda_n < \infty$. Note that the matrix $\vec{a}$ has the unit of a diffusion (i.e. $\text{m}^2 \cdot \text{s}^{-1}$). The basis functions can be constructed by numerical or observational data and various forms of noise can be used in the simulations, see e.g. \cite{Brecht2021b}. This allows us to formulate the RSW-LU as:
\begin{align}
    \df_t \vec{u} & = \overbrace{\Big( - \vec{u} \adv \vec{u} - \vec{f} \times \vec{u} - g \grad h \Big)}^{:=\operatorname{\bf det}}\, \dt 
+ \overbrace{\Big( \alf \div \div (\mathbf{a} \mathbf{u})\, \dt - \sdbt \adv \mathbf{u} \Big)}^{:=\operatorname{\bf sto}^V },
\label{eq:rsw_determ}\\
	\df_t h       & = - \div (\vec{u} h)\, \dt 
	+ \underbrace{\Big( \alf \div \div (\vec{a} h)\, \dt - \sdbt \adv h \Big)}_{:=\operatorname{\bf sto}^h },\label{eq:rsw_determ1}
\end{align}
where $f$ is the Coriolis parameter and $g$ the gravitational acceleration
and where $\alf \div \div (\mathbf{a} \mathbf{\theta}) = \alf (\div \mathbf{a})\cdot \nabla \theta + \alf \div (\mathbf{a} \nabla \theta)$
with $\div \mathbf{a} = \sum_i \frac{\partial a_{ik}}{\partial x_i}$ for the coordinate functions $\theta = u, v, h$.

As suggested in \cite{Brecht2021b}, our approach is to discretize the deterministic parts (the $\operatorname{det}$ and $-\nabla\cdot (\mathbf{u}h)$ terms) with a variational discretization \cite{BauerRSW2019}
and the stochastic terms with standard finite difference operators. Then, we add an energy preserving Casimir dissipation \cite{Brecht2021a} or biharmonic Laplacian dissipation.

\section{Discretization of RSW-LU with Casimir dissipation}
\label{sec_discretization}

In this section we first introduce the variational principle that underpins
the derivation of the discrete dissipative variational equations before we combine this 
deterministic scheme with approximations of the stochastic terms such that 
together they form a semi-discrete approximation of the RSW-LU equations.

\subsection{Discrete variational equations}

\begin{figure}
\centering
\begin{minipage}{6cm}
\includegraphics[scale=0.55,angle=0]{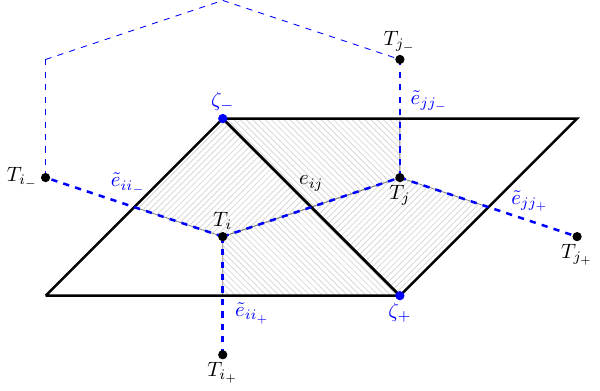}
\end{minipage}
\begin{minipage}{5cm}
{\renewcommand{\arraystretch}{1.2}
\begin{tabular}{|c|l|}
\hline
$T_i$        & triangle with area $\Omega_{ii}$        \\ \hline
$\zeta_\pm$  & dual cell       \\ \hline
$e_{ij}$ & primal edge $T_i\cap T_j$   \\ \hline
$\tilde{e}_{ij}$ &dual edge $\zeta_+\cap \zeta_-$ \\ \hline
$\mathbf{n}_{ij}$ & $e_{ij}$ normal vector towards $T_j$ \\ \hline
$h_i$ & waterdepth on $T_i$ \\ \hline
$\bar{h}_{ij}$ & $\frac{1}{2} (h_i+h_j)$  \\ \hline
$V_{ij}$ & $(\mathbf{u}\cdot \mathbf{n})(e_{ij})$ \\ \hline 
\end{tabular}
}
\end{minipage}
\caption{Notation for 2D simplicial mesh.} \label{fig:notation}
\vspace{-1em}
\end{figure}

We consider a two-dimensional simplicial mesh $\M$  with $n$ cells on the fluids domain, where  triangles ($T$) are used as the primal grid, and the circumcenter dual ($\zeta$) as the dual grid. 
The notation is explained in Fig. \ref{fig:notation}.
\vspace{-1em}

Following \cite{BauerRSW2019,Brecht2021b}, let matrix $A$ approximate the vector field $\vec{u}$ and vector
$h\in \R^n$ the fluid depths from RSW-LU. Here, piecewise constant functions on the domain $\M$ are represented by vectors $F\in\R^n$, with value $F_i$ on cell $i$ being the cell average of the continuous function on cell $i$. The space of discrete functions is denoted by $\Omega^0_d(\M)$, and the space of discrete densities $\Den_d(\M)\simeq \R^n$ is defined as the dual space to $\Omega^0_d(\M)$ relative to the pairing: $\langle F, G\rangle_0 = F^\top \Omega G$, 
where $\Omega$ is a nondegenerate diagonal $n \times n$ matrix with elements $\Omega_{ii}=\mathrm{Vol}(T_i)$, the area of triangle $T_i$.

$A$ is an element of the Lie algebra $\mathfrak{d}(\M)=\{ A\in \mathfrak{gl}(n) \mid A\cdot \mathbf{1}=0\}$ with the matrix commutator $[A,B]=AB-BA$ as the Lie bracket, where $\mathfrak{gl}(n)$ is the Lie algebra of $n\times n$ real matrices. To obtain 
approximations of vector fields we consider a subspace  $ \mathcal{R} \subset \mathfrak{d}(\M)$ in which 
$\mathcal{R} = R_1 \cap R_2$, $R_1 = \big\{ A \in \mathfrak{d}(\M)\mid A^\top \Omega + \Omega A^\top \text{ is diagonal }\big\}$ and $R_2=\big\{ A \in \mathfrak{d}(\M)\mid A_{ij}=0 ~~~\forall j\notin N(i)\big\}$
with $N(i)$ being the set of cells sharing an edge with the cell $T_i$. 
Further, $A_{ij}= -\frac{|e_{ij}|}{2\Omega_{ii}}  V_{ij}, \; j\in N(i)$ and $ A_{ii} =  \frac{1}{2\Omega_{ii}}\sum_{k\in N(i)} |e_{ik}|V_{ik}$, where the former is an off diagonal element associated to the triangular edge with normal pointing from cells $i$ to $j$, while the latter is a diagonal element associated to cell $i$. We identify the dual space $\mathcal{R}^*$ with the space $\Omega^1_d(\M)$ of discrete one-forms relative to the duality pairing on $\mathfrak{gl}(n)$: 
$\langle L, A\rangle_1 = \text{Tr}(L^\top \Omega A)$. To obtain an element in $\mathcal{R}^*$, 
$P\colon \mathfrak{gl}(n)\to \Omega^1_d(\M)$ is used, defined by
$P(L)_{ij}=\frac{1}{2}(L_{ij}-L_{ji}-L_{ii}+L_{jj})$, 
which satisfies $ \left\langle L, A \right\rangle _1 = \left\langle P(L), A \right\rangle _1$, for all $A \in \mathcal{R}$, cf. \cite{BauerRSW2019} for details.

Then, let $\ell\colon  \mathfrak{d}(\M) \times \Den_d(\M)\to \R$ be a semi-discrete Lagrangian and $C\colon  \mathfrak{d}(\M)  \times \Den_d(\M)\to \R$ be a semi-discretized approximation of a Casimir.
\begin{theorem}[\textbf{Discrete dissipative variational equations}]\label{th:discreteVarEq}
	For a semi-discrete Lagrangian $\ell(A,h)$, the curves $A(t),h(t) \in \mathcal{R}$ are critical for the variational principle of Eq. (28) in \cite{Brecht2021a} if and only if they satisfy
	\begin{equation}	\label{eq:projectionCas}
		P \left(	\frac{d}{dt}\frac{\delta \ell}{\delta A}+ \mathcal{L}_A\Big(\frac{\delta \ell}{\delta A}\Big)- \theta\mathcal{L}_A\Big(h  {\mybrL}\frac{\delta C}{\delta M},A {\mybrR}^\flat\Big)+ h\frac{\delta \ell}{\delta h}^\top\right)_{ij}=0,
	\end{equation}				
	where $\mathcal{L}$ is the discrete analog to the Lie derivative $\mathfrak{L}$ and it is defined by the commutator via the following relation 
$			\left\langle \mathcal{L}  _A M, B \right\rangle _1 = \left\langle M, [A,B] \right\rangle _1$. 
The proof can be found in \cite{Brecht2021a}.
\end{theorem}
Note that for $A,B\in \mathcal{R}$ we have $[A,B]_{ij}=0$ for all $j \in N(i)$. Since elements of $\mathcal{R}$ are zero for non-neighboring cells, we get $[\mathcal{R},\mathcal{R}]\cap \mathcal{R}=\{0\}$. In particular $[ \mathcal{R} , \mathcal{R} ] \neq \mathcal{R} $ hence the subspace $ \mathcal{R}  \subset \mathfrak{d}(\M)$ corresponds to a nonholonomic constraint. Consequently, we need to define a discrete commutator $\llbracket , \rrbracket$ such that $\llbracket A, B \rrbracket \in \mathcal{R}$, so that we can directly apply the definition of the discrete flat operator to the discrete commutator, cf. \cite{Brecht2021b} for details on this construction.

In the case of the RSW equations, we have the discrete Lagrangian
\begin{align}\label{eq:discrteLagrangian}
	\begin{split}
		\ell(A, h) &= \frac{1}{2}\sum_{i,j=1}^n h_i A^\flat_{ij}A_{ij} \Omega_{ii} 
		+\sum_{i,j=1}^n h_i R^\flat_{ij} A_{ij} \Omega_{ii}
		-\frac{1}{2}\sum_{i=1}^n g(h_i+(\eta_b)_i)^2\Omega_{ii}, 
\end{split}
\end{align} 
where $R$ is the vector potential of the angular velocity of the Earth and 
$\flat$ is the discrete flat operator \cite{BauerRSW2019}. The functional derivatives are given by $\frac{\delta \ell}{\delta A}_{ij} = h_i (A_{ij}^\flat+R_{ij}^\flat)$
and $\frac{\delta \ell}{\delta h}_i =\frac{1}{2}\sum_j A^\flat_{ij}A_{ij}+\sum_j R^\flat_{ij}A_{ij}-g(h_i+(\eta_b)_i)$.
\smallskip

We define the discrete enstrophy Casimir (i.e. the discrete potential enstrophy) as 
\begin{equation}
\mathcal{C}(M,h) =\frac{1}{2} \sum_{\zeta } h_{\zeta} \Big(q(M,h)_\zeta\Big)^2 |\zeta| ,
\end{equation}
for $M=\frac{\delta \ell}{\delta A}$, with discrete potential vorticity $q(M,h)_\zeta  =\frac{(\text{Curl}~V)_\zeta + f}{h_\zeta}$ 
and discrete fluid depth (mass) $h_\zeta  = \sum_{T_i \cap \zeta \neq \emptyset} \frac{|T_i\cap \zeta|}{|\zeta|} h_i$,
 which are discrete realizations of the mass, potential vorticity and 
potential enstrophy Casimirs of the continuous RSW equations. The functional derivative of $\C$ reads 
$	\left( \frac{\delta \C}{\delta M}\right)_{ij}
=-\frac{|e_{ij}|}{2\Omega_{ii}} \frac{2~\text{Grad}_t ~q}{\overline{h}_{ij}}$ with $\text{Grad}_t$ defined below, 
 see computations in \cite{Brecht2021a}.

\subsection{Semi-discrete RSW-LU scheme}

The resulting system of semi-discrete equations can be split into a \textbf{det}erministic, \textbf{diff}usive and \textbf{sto}chastic part. In the momentum equation we include the two different diffusion methods (i.e. Casimir or biharmonic Laplacian) which can be applied by setting $\theta>0$ or $\nu>0$:
\begin{align}
\begin{split}\label{eq:detdiffsto}
    \partial_t V_{ij} &= - \operatorname{\bf det}_{ij} - \theta \operatorname{\bf diff}_{ij}^{\text{CD}} - \nu \operatorname{\bf diff}_{ij}^{\text{BD}} + \operatorname{\bf sto}_{ij}^V,
    \\
    \partial_t h_{i}&=-\operatorname{Div}(\bar{h}V)_i + \operatorname{\bf sto}_i^h .
\end{split}\end{align}
Note that the terms $\operatorname{\bf det}_{ij}$ and $\theta \operatorname{\bf diff}_{ij}^{\text{CD}}$ in the momentum equation together with $\operatorname{Div}(\bar{h}V)_i$ in the continuity equation constitute the Casimir-dissipation stabilized RSW equations that result from Eqn.~\eqref{eq:projectionCas}.

Before specifying these terms, we define the following finite differences operators:
\begin{align*}
(\operatorname{Grad}_n F)_{ij} &= \frac{F_{T_j}-F_{T_i}}{|\tilde{e}_{ij}|} ,
&&
(\operatorname{Div} V)_{i}=\frac{1}{|T_i|}\sum_{k\in\{j,i_-,i_+\}} |e_{ik}| V_{ik}, 
\\
(\operatorname{Grad}_t F)_{ij} &= \frac{F_{\zeta_-}-F_{\zeta_+}}{|e_{ij}|}, 
&&
(\operatorname{Curl} V)_{\zeta}=\frac{1}{|\zeta|}\sum_{\title{e}_{nm}\in \partial \zeta} |\title{e}_{nm}| V_{nm}.
\end{align*}
Predicting only the normal velocity component $V_{ij}$ for all edges, for some terms we need the full velocity field; for this we apply the reconstruction: \\
$\mathbf{u}_i = \frac{1}{\Omega_{ii}} \sum_{k\in\{j,i_-,i_+\}}|e_{ik}|(\mathbf{x}_{e_{ik}}-\mathbf{x}_{T_i})V_{ik}$.

The various deterministic terms in Eq. \eqref{eq:detdiffsto} are given by: 
\begin{small}
\begin{align*}
	 \operatorname{\bf det}_{ij} := &
   -  \frac{(\text{Curl } V + f)_{\zeta_{-}} }{\overline{h}_{ij} |\tilde{e}_{ij}| } 
	\Big(   \frac{|\zeta_{-} \cap T_i|}{2 \Omega_{ii}  } \overline{h}_{ji_-} |e_{ii_-}|   V_{ii_-} 
	+  \frac{|\zeta_{-}  \cap T_j|}{2 \Omega_{jj}  } \overline{h}_{ij_-}  |e_{jj_-}|   V_{jj_-} \Big) \\
	& +  \frac{ (\text{Curl } V + f)_{\zeta_{+}} }{\overline{h}_{ij}  |\tilde{e}_{ij}| } 
	\Big(  \frac{|\zeta_{+}  \cap T_i|}{2 \Omega_{ii} }\overline{h}_{ji_+} |e_{ii_+}|   V_{ii_+} + \frac{|\zeta_{+}  \cap T_j|}{2 \Omega_{jj}   } \overline{h}_{ij_+}   |e_{jj_+}|   V_{jj_+} \Big)
\\	
 & - \frac{1}{2}\Big(\text{Grad}_n \Big(\sum_{k\in \{j, i_-, i_+\} }\frac{|\tilde{e}_{ik}|~|e_{ik}|(V_{ik})^2}{2\Omega_{kk}}\Big)\Big)_{ij}
-	 g (\text{Grad}_n~(h+\eta_b))_{ij} \ ,
\\
 \operatorname{\bf diff}_{ij}^{\text{CD}} := &
 -  \frac{(\text{Curl } \widetilde{W})_{\zeta_{-}} }{\overline{h}_{ij} |\tilde{e}_{ij}| } 
	\Big(   \frac{|\zeta_{-} \cap T_i|}{2 \Omega_{ii}  } \overline{h}_{ji_-} |e_{ii_-}|   V_{ii_-} 
	+  \frac{|\zeta_{-}  \cap T_j|}{2 \Omega_{jj}  } \overline{h}_{ij_-}  |e_{jj_-}|   V_{jj_-} \Big) \\
	& +  \frac{ (\text{Curl }\widetilde{W})_{\zeta_{+}} }{\overline{h}_{ij}  |\tilde{e}_{ij}| } 
	\Big(  \frac{|\zeta_{+}  \cap T_i|}{2 \Omega_{ii} }\overline{h}_{ji_+} |e_{ii_+}|   V_{ii_+} + \frac{|\zeta_{+}  \cap T_j|}{2 \Omega_{jj}   } \overline{h}_{ij_+}   |e_{jj_+}|   V_{jj_+} \Big) 
 \\
      & +\frac{(2 \widetilde{W}_{ij})}{\overline{h}_{ij}  } 
	      \frac{\Div(V\overline{h})_i+\Div(V\overline{h})_j}{2}   
 - \frac{1}{2}\Big(\text{Grad}_n \Big( \!\!\!\!\! \sum_{k\in \{j, i_-, i_+\} }  \!\!\!\!\! \frac{|\tilde{e}_{ik}|~|e_{ik}|(V_{ik}\widetilde{W}_{ik})}{\Omega_{kk}}\Big)\Big)_{ij} \ ,
  \\
  \operatorname{\bf diff}_{ij}^{\text{BD}} := & \ \text{Lap}(\text{Lap}(V))_{ij}  \ \text{with} \ 
   \text{Lap}(V)_{ij} := \left(\text{Grad}_n (\operatorname{Div}(V))_{ij} - \text{Grad}_t (\text{Curl} V)_{ij}  \right)  ,
\end{align*}
\end{small}
\vspace{-1em}

\noindent
where $\widetilde{W}_{ij}$ is given by
\begin{align*}
    \widetilde{W}_{ij} =&\left(\frac{\delta C}{\delta M}\right)_{ij}\Big(\frac{\Div(V)_i+\Div(V)_j}{2}\Big)
	  -V_{ij}\Big(\frac{\Div\left(\frac{\delta C}{\delta M}\right)_i+\Div\left(\frac{\delta C}{\delta M}\right)_j}{2}\Big) \\
    &~~~ -\text{Grad}_t~\Big((\frac{\delta C}{\delta \mathbf{m}}_{\zeta}\times \mathbf{u}_\zeta)\cdot \mathbf{k}_\zeta\Big)_{ij} \  \text{ with } \ \mathbf{u}_\zeta = \sum_{i\in N(\zeta)} \frac{|\zeta\cap T_i|}{|\zeta|} \mathbf{u}_i. 
\end{align*}

The stochastic term $\operatorname{\bf sto}_{ij}^V := \operatorname{\bf sto}^V (ij) \cdot \mathbf{n}_{ij}$ follows from an evaluation of the stochastic terms in \eqref{eq:rsw_determ} at the edge midpoints $ij$ and a projection onto the edges' normals, and $\operatorname{\bf sto}_{i}^h := \operatorname{\bf sto}^h (i)$ from an evaluation of the terms in \eqref{eq:rsw_determ1} at cell centers $i$. They further apply a reconstruction of the 3D velocity field in Cartesian coordinates and the evaluation of the partial derivatives using the approximation:
$ (\partial_{x_m} F)_{ij} = (\operatorname{Grad}_n F)n^m_{ij}+(\operatorname{Grad}_t F)t^m_{ij}, \ m=1,2,3$, where $n$ and $t$ are the edge normal and tangential vector components. Full details 
can be found in \cite{Brecht2021b}.

\section{Numerical results}
\label{sec_numerics}

We illustrate on a selected test case (TC), namely a barotropically unstable jet at mid-latitude on the sphere, proposed by \cite{Galewsky2004} and here referred to as Galewsky TC, the superior performance of our Casimir-dissipation (\textbf{CD}) stabilized stochastic flow model over a comparable version that applies instead a conventional biharmonic Laplacian diffusion (\textbf{BD}). This test case is particularly suited to study the performance of stochastic schemes for low resolution (LR) grids: as small scale features trigger and determine the large scale flow pattern, the former have to be either well resolved (normally with sufficiently high resolution) or accurately represented by the stochastic model. Here, comparable deterministic LR simulations usually fail to show the correct large scale pattern, cf. \cite{Brecht2021b}. Therein, also the initialization of the water depths and velocity fields can be found.

As this manuscript compares the use of two different dissipation methods, one might ask if there is the need for dissipation at all. This is a valid question, especially when considering that we use a spatially energy conserving discretization of a consistent stochastic flow model (see above) which should, in principle, be stable even without diffusion. 
For the given test case and integrations up to 20 days, our RSW-LU scheme is indeed numerically stable, but we did not explore longer integration times because some kind of dissipation is in general anyway needed. This is mostly due to the accumulation of small scale potential enstrophy (PE) at the grid level -- caused in 2D simulations by the PE cascade towards small scales. If not dissipated, either by CD or BD diffusion, this accumulation of PE leads to noisy and nonphysical fields or even to numerical instabilities (shown in \cite{Brecht2021a} for deterministic models, but also valid for stochastic ones).

\begin{figure}[h]\centering
\begin{tabular}{cc}
 \hspace{-1em}\includegraphics[scale = 0.24]{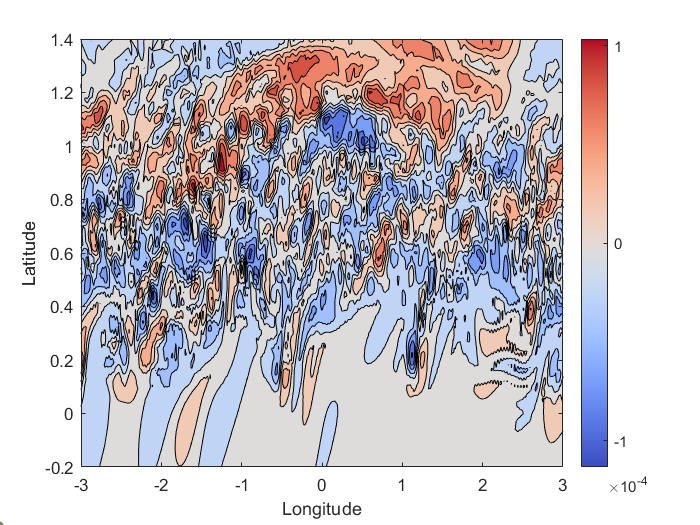} &
 \hspace{-1em} \includegraphics[scale = 0.29]{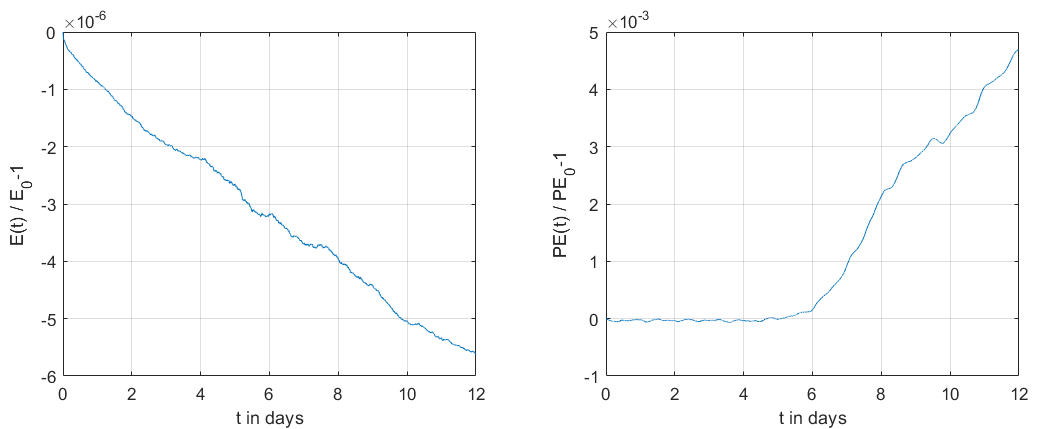} 
\end{tabular}
\caption{PV field at day 12 (left), and total energy (middle) and potential enstropy (right) development up to day 12 for the stochastic RSW-LU without CD or BD dissipation. The contour interval in the left panel is $2 \times 10^{-5} \rm{s^{-1}}$.}
\label{fig_nocdnobd}
\end{figure}

To support these points, we first show in Figure~\ref{fig_nocdnobd} the results of a simulation over 12 days of the Galewsky TC on a grid with 20480 triangles and without CD or BD dissipation.
Looking at the left panel, we see that the potential vorticity (PV) field exhibits mostly nonphysical small scale noise that results from the accumulation of PE at the grid scale and does not agree with the physics of fluid flows. Looking at the total PE time evolution (right), we observe that with the beginning of day 5, when the large scale flow pattern starts to emerge, PE starts to increase significantly. This increase might eventually lead to a crash of the numerical simulation. Using a spatially energy conserving scheme, the slight drift in total energy (middle) is due to the time integration scheme of our stochastic system. For further discussions and more details about how total energy (E) and PE are computed, see \cite{Brecht2021b}.

\begin{figure}[h]\centering
\begin{tabular}{c}
 \includegraphics[scale = 0.35]{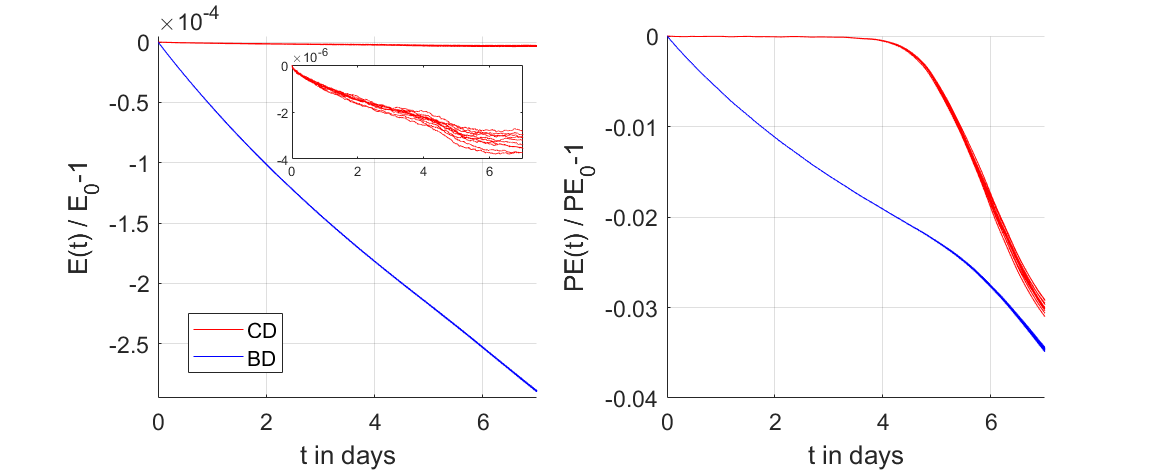}
\end{tabular}
\caption{Time evolution of total energy (left) and potential enstrophy (right) for an integration time of 7 days for CD (red) and BD (blue) for all 20 ensemble members.}
\label{fig_en_enst}
\end{figure}

Having illustrated the necessity of dissipating PE, we next compare the BD and CD dissipation methods. To make this comparison fair, we proceed as follows. We first tune the BD diffusion coefficient: we let the low resolution BD simulation run for 12 days with a sufficiently low value for $\nu$ and then increase it until we suppress numerical ringing in the PV field. This gives us an optimal value of about $\nu = 3.1 \times 10^{16}\rm{m^4/s}$ for this TC. As can be inferred from Figure~\ref{fig_en_enst}, this dissipates beside PE also the total energy E. Considering the frequency spectrum for PE (calculated using the method from \cite{Brecht2021a} but not shown here), we now choose the Casimir dissipation coefficient $\theta =  5\cdot 10^{21}\ \rm{m^5 s}$ such that both PE frequency spectra agree, especially at small scales, which indicates that both schemes sufficiently well dissipate the small scale PE noise.

\begin{figure}
\begin{tabular}{c}
\hspace*{-2em} \includegraphics[scale = 0.40]{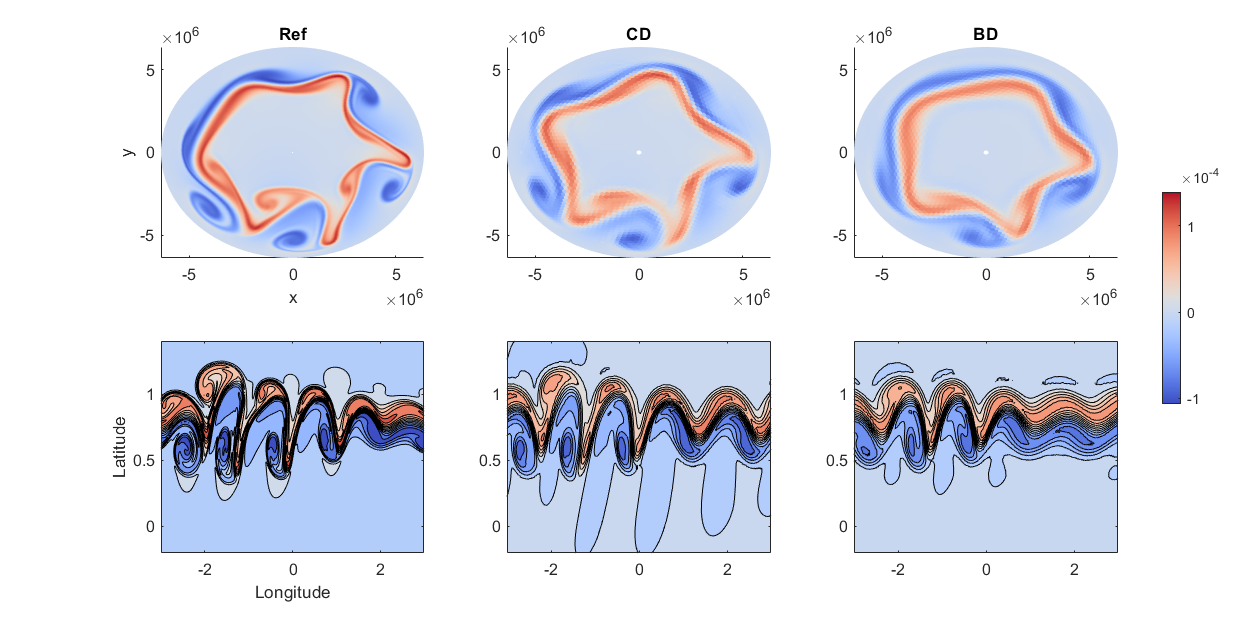}  
\end{tabular}
\caption{Ensemble mean of PV fields at day 6 for high resolution (Ref) simulation with 327680 triangles (first column) and for low resolution CD (second column) and BD (third column) simulations. Upper row shows the northern hemisphere from above, lower row shows the projections of these fields onto latitude-longitude grid representations. The contour interval in the lower row pictures is $2 \times 10^{-5} \rm{s^{-1}}$.  
}
\label{fig_no_cd_bd}
\end{figure}

In Figure~\ref{fig_no_cd_bd}, we compare for the Galewsky TC the PV fields at day 6 obtained from a high resolution reference large-eddy simulation (Ref) on a grid with 327680 triangles (as used in \cite{Brecht2021b}) (first column) with two low resolution runs where one uses only CD dissipation (second column) and the other only BD dissipation (third column). Note that these fields show the ensemble mean of an ensemble run with 20 members. 
We observe that the mean flow field of the CD ensemble is closer to Ref than that of the BD ensemble: in the former case, the eddies have been better developed while the distances between vortex (blue) and anti-vortex (red) pairs are larger, hence better agree with those from the reference run. 
In contrast, vortex anti-vortex pairs are less well developed in the BD runs. As the latter dissipate, besides PE, also total energy, it seems that using BD dissipation removes also dynamics from the system. 

This conclusion is supported by the time evolution of E and PE of the 20 ensemble members, shown in Figure~\ref{fig_en_enst}. For both quantities, we show in red results for CD and in blue for BD. All ensemble members of the CD runs preserve total energy at the order of $10^{-6}$. Comparing them with the corresponding curve in Fig.~\ref{fig_nocdnobd} confirms that our implementation of CD does not diffuse energy (diffusion is only due to time integration). Besides, all members of the CD ensemble diffuse almost the same amount of PE. In contrast, all members of the BD runs diffuse besides PE also a substantial amount of energy (at the order of $10^{-4}$), leading to the reduces dynamics as stated above.

\section{Conclusions}
\label{sec_conclusion}

We presented a novel structure preserving discretization of a Location Uncertainty version of stochastic rotating shallow water equations, applying an energy conserving Casimir diffusion (CD) to dissipate small scale potential enstrophy. 
We compared this scheme with a version that applies biharmonic Laplacian diffusion (BD) for stabilization, which is usually used for such models. Our results show clearly that the CD runs achieve a spatial conservation of energy, which, in contrast, is strongly dissipated for BD. Consequently, compared to BD, the dynamics of CD simulations are significantly closer to a high resolutions reference solution, underpinning the use of our CD over standard BD diffusion.

\subsubsection{Acknowledgements} 
RB is funded by the Deutsche Forschungsgemeinschaft (DFG, German Research Foundation) – Project-ID 274762653 – TRR 181.

\end{document}